\documentclass[final]{jotart}

\usepackage[mathscr]{eucal}
\usepackage{amsmath}
\usepackage{amssymb}

\theoremstyle{proclaim}
\newtheorem{theorem}{Theorem}[section]
\newtheorem*{theorem*}{Theorem}
\newtheorem{lemma}[theorem]{Lemma}
\newtheorem{corollary}[theorem]{Corollary}
\newtheorem{proposition}[theorem]{Proposition}
\theoremstyle{fancyproclaim}

\theoremstyle{statement}
\newtheorem{remark}[theorem]{Remark}
\newtheorem{definition}[theorem]{Definition}
\newtheorem{example}[theorem]{Example}
\theoremstyle{fancystatement}

\numberwithin{equation}{section}

\providecommand{\AMS}{$\mathcal{A}$\kern-.1667em%
\lower.25em\hbox{$\mathcal{M}$}\kern-.125em$\mathcal{S}$}

\newcommand{\be}{\begin{equation}\label}
\newcommand{\ee}{\end{equation}}
\newcommand{\bq}{\begin{equation*}}
\newcommand{\eq}{\end{equation*}}
\newcommand{\ba}{\begin{align*}}
\newcommand{\ea}{\end{align*}}
\newcommand{\bp}{\begin{proof}}
\newcommand{\ep}{\end{proof}}
\newcommand{\bL}{\begin{lemma}\label}
\newcommand{\eL}{\end{lemma}}
\newcommand{\bP}{\begin{proposition}\label}
\newcommand{\eP}{\end{proposition}}
\newcommand{\bC}{\begin{corollary}\label}
\newcommand{\eC}{\end{corollary}}
\newcommand{\bT}{\begin{theorem}\label}
\newcommand{\eT}{\end{theorem}}
\newcommand{\bR}{\begin{remark}\label}
\newcommand{\eR}{\end{remark}}
\newcommand{\bD}{\begin{definition}\label}
\newcommand{\eD}{\end{definition}}

\newcommand{\bE}{\begin{example}\label}
\newcommand{\eE}{\end{example}}
\newcommand{\bQ}{\begin{question}\label}
\newcommand{\eQ}{\end{question}}

\DeclareMathOperator{\diag}{diag}

\begin{document}

\commby{K.Davidson}
\date{September 02, 2011}
\revision{February 06, 2012}
\title[Subideals of Operators]{SUBIDEALS OF OPERATORS}

\author{SASMITA PATNAIK \protect AND GARY WEISS}
\address{Sasmita Patnaik, Department of Mathematics, University of Cincinnati, OH,45221-0025 , USA}
\email{sasmita\_19@yahoo.co.in}
\address{Gary Weiss, Department of Mathematics, University of Cincinnati, OH, 45221-0025, USA}
\email{gary.weiss@math.uc.edu}

\begin{abstract}
A subideal is an ideal of an ideal of $B(H)$ and a principal subideal is a principal ideal of an ideal of $B(H)$.
We determine necessary and sufficient conditions for a principal subideal to be an ideal of $B(H)$.
This generalizes to arbitrary ideals the 1983 work of Fong and Radjavi characterizing principal subideals of the ideal of compact operators that are also ideals of $B(H)$.
We then characterize all principal subideals.
We also investigate the lattice structure of subideals as part of the general study of ideal lattices such as the often studied lattice structure of ideals of $B(H)$.
This study of subideals and the study of elementary operators with coefficient constraints are closely related.
\end{abstract}



\maketitle

\section{INTRODUCTION}\label{S: 1}

This paper investigates the subideal structure of $B(H)$ following the spirit of Calkin's well-known singular number characterization of ideals of $B(H)$ (i.e., henceforth called $B(H)$-ideals) \cite{C41}.
A subideal is an ideal of an ideal $J$ (henceforth called a $J$-ideal) for the $B(H)$-ideal $J$.
Recall for general rings, an ideal is an additive commutative subgroup which is closed under left and right multiplication by elements of the ring.
Ideals in the ring $B(H)$ are ubiquitous throughout operator theory.
Some well-known $B(H)$-ideals are the compact operators $K(H)$, the finite rank operators $F(H)$, principal ideals, Banach ideals, the Hilbert-Schmidt class, the trace class, Orlicz ideals, Marcinkiewicz ideals and Lorentz ideals, to name a few.
Definitions of these ideals may be found in \cite{DFWW}.
Recall also that all proper $B(H)$-ideals lie in $K(H)$.
Here and throughout this paper $H$ denotes a separable infinite-dimensional complex Hilbert space, $B(H)$ the algebra of all bounded linear operators on $H$, and $\mathbb {C, R, N, Z}$, respectively, the classes of complex numbers, real numbers, positive integers and integers.\\

There are three natural kinds of \emph{principal $J$-ideals}, namely, the classical principal $J$-ideals $(S)_J$ which we call principal linear $J$-ideals; principal $J$-ideals $\left<S\right>_J$ and principal real linear $J$-ideals $(S)_J^{\mathbb{R}}$ (Definition \ref{D:d1}).
Standard notation then dictates that we denote $(S) = (S)_{B(H)}$.
\emph{It is immediate that $B(H)$-ideals are always $J$-ideals but, as we shall see later, often not conversely (\cite{FR83}, see also Example} \ref{E:1.1} below).\\

The main results of this paper,
generalizing the 1983 work of Fong and Radjavi \cite{FR83} and characterizing the principal subideals of $B(H)$, are summarized in the following two theorems.

For compact operators $S, T$, $s(S)$ denotes the sequence of singular numbers (s-numbers) for $S$,
and the product $s(S)s(T)$ denote their pointwise product.

\bT{T:1}
For $S \in J$, the following are equivalent.\\

\nr{i} Any of the three types of principal $J$-ideals generated by $S$, $(S)_J$, $\left<S\right>_J$ or $(S)_J^{\mathbb{R}}$,

~is a $B(H)$-ideal.

\nr{ii} The principal $B(H)$-ideal
$(S)$ is $J$-soft, i.e., $(S)$ = $J(S)$ (equivalently, $(S)$ = $(S)J$).

\nr{iii} $S = AS + SB + \displaystyle{\sum_{i=1}^{m}}A_{i}SB_{i}$ for some $A,\,B,\,A_{i},\, B_{i} \in J,~ m \in \mathbb{N}$.

\nr{iv} $s(S) = \textrm{O}(D_{k}(s(S))s(T))$ for some $T \in J$ and $ k \in \mathbb{N}$.\\
\eT

Denoting $JS + SJ + J(S)J := \{AS + SB + CS'D \mid A, B, C, D \in J, S' \in (S)\},$
\bT{T:2}
The principal $J$-ideal, the principal linear $J$-ideal and the principal real linear $J$-ideal generated by $S \in J$ are respectively given by
\begin{equation*}
\left<S\right>_J = \mathbb Z S + JS + SJ + J(S)J
\end{equation*}
\begin{equation*}
(S)_J = \mathbb C S + JS + SJ + J(S)J
\end{equation*}
\begin{equation*}
(S)_J^\mathbb R = \mathbb R S + JS + SJ + J(S)J.
\end{equation*}
\begin{center}
\begin{equation*}
\text{So} \quad J(S)J \subseteq JS + SJ + J(S)J \subseteq \left<S\right>_J \subseteq (S)_J^\mathbb R \subseteq (S)_J \subseteq (S)
\end{equation*}
\end{center}
which first two, $J(S)J$ and $JS + SJ + J(S)J$ respectively, are a $B(H)$-ideal and a $J$-ideal.
Consequently, each of these three kinds of principal $J$-ideals have the common $B(H)$-ideal ``nucleus'' $J(S)J$,
with the common $J$-ideal $JS + SJ + J(S)J$ containing it.

All these principal $J$-ideals, $\left<S\right>_J \subseteq (S)_J^\mathbb R  \subseteq (S)_J $, are distinct except under the following equivalent conditions.
They all collapse to merely
\begin{equation*}
J(S)J = (S) = \left<S\right>_J = (S)_J = (S)_J^\mathbb R
\end{equation*}
if and only if the principal $B(H)$-ideal $(S)$ is $J$-soft \\
(that is, $(S) = J(S)$ (Definition \ref{D:1}) in which case $(S) = J(S)= J(S)J$)\\
if and only if any, and hence all of them, is a $B(H)$-ideal.
\eT
\emph{Background}

In 1941, Calkin \cite{C41} characterized $B(H)$-ideals via his lattice preserving isomorphism between $B(H)$-ideals and characteristic sets $\Sigma \subseteq c_{0}^*$:
$I \rightarrow \Sigma(I)$ induced by $I \owns X \rightarrow s(X) \in \Sigma(I)$.
Here $c_{0}^{*}$ denotes the cone of non-negative sequences decreasing to zero; characteristic sets $\Sigma$ are those subsets of $c_{0}^{*}$ that are additive, hereditary (solid) and ampliation invariant (invariant under each $m$-fold ampliation $D_m\xi :=\,<\xi_{1},\cdots, \xi_{1}, \xi_{2}, \cdots, \xi_{2}, \cdots>$ with each entry $\xi_{i}$ repeated $m$ times);
the characteristic set $\Sigma(I) := \{\eta \in   c_{0}^* \mid \diag \eta \in I\}$;
and $s(X)$ denotes the $c_{0}^*$-sequence of s-numbers of compact operator X.\\

Motivated by this characterization, a natural question to ask and the subject of this paper is:
\begin{center}
\emph{What can be said about subideals, i.e., is it possible to characterize them in some way?}
\end{center}

In 1983, Fong and Radjavi \cite{FR83} investigated those subideals that are principal linear $K(H)$-ideals, perhaps in part because of the distinguished role $K(H)$ plays as the unique norm closed proper $B(H)$-ideal.
Though unstated there, ``linearity'' was assumed as recently clarified to us (private communications).
They found principal linear $K(H)$-ideals that are not $B(H)$-ideals (Example \ref{E:1.1})
by determining necessary and sufficient conditions for a  principal linear $K(H)$-ideal to be a $B(H)$-ideal
\cite[Theorem 2]{FR83}. And in doing so, at least for the authors of this paper, they initiated the study of subideals.

\begin{theorem*} \cite[Theorem 2]{FR83} Let $T$ be a compact operator of infinite rank and let $P = (T^*T)^{\frac{1}{2}}$. Let $\mathcal{T}$ and $\mathcal{P}$ be the ideals in $K(H)$ generated by $T$ and $P$, respectively. Then the following are mutually equivalent.

\nr{i} $\mathcal{T}$ is an ideal in $B(H)$.

\nr{ii} $\mathcal{P}$ is an ideal in $B(H)$.

\nr{iii} $T =   A_1TB_1 + \dots + A_kTB_k$ for some $k$ and some $A_i \in K(H)$, $B_i \in B(H)$.

\nr{iv} $T =   A_1TB_1 + \dots + A_kTB_k$ for some $k$ and some $A_i \in K(H)$, $B_i \in K(H)$.
\end{theorem*}

Fong and Radjavi proved this via the positive case employing the Lie ideal condition (ii) below.

\begin{theorem*}\cite[Theorem 1]{FR83} Let $P$ be a positive compact operator of infinite rank, and let $\mathcal{I}$ be the ideal in $K(H)$ generated by $P$. Then the following are equivalent.\\
\nr{i} $\mathcal{I}$ is an ideal in $B(H)$.

\nr{ii} $\mathcal{I}$ is a Lie ideal in $B(H)$.

\nr{iii} $P = A_1PA_1^{*} + \dots + A_kPA_k^{*}$ for some $k$ and some compact operators $A_i$.

\nr{iv} $P = A_1PB_1 + \dots + A_kPB_k$ for some $k$ and some compact operators $A_i$ and $B_i$.

\nr{v} $P = A_1PB_1 + \dots + A_kPB_k$ for some $k$, \\
\indent where $A_i$, $B_i \in B(H)$ and either the $A_i$ or the $B_i$ are compact.

\nr{vi} For some integer $k > 1$, $s_{nk}(P)= o(s_n(P))$ as $n \rightarrow \infty$.
\end{theorem*}

\bE{E:1.1}
Condition (vi) of \cite[Theorem 1]{FR83} shows that if the singular number sequence of the operator $P$ is given by $s(P)= \left<\frac{1}{2^n}\right>$,  then the principal linear $K(H)$-ideal generated by $P$ is a $B(H)$-ideal.
But if $s(P)= \left<\frac{1}{n}\right>$, then the principal linear $K(H)$-ideal generated by $P$ is not a $B(H)$-ideal.
\eE

\emph{In summary}

This paper fully generalizes \cite[Theorem 2]{FR83} from principal $K(H)$-ideals to arbitrary principal $J$-ideals and all but the Lie ideal condition in \cite[Theorem 1]{FR83}.
We investigate all three types of principal $J$-ideals, whereas Fong-Radjavi considered only the principal linear $J$-ideals and for only the case $J=K(H)$ \cite{FR83}.
We determine necessary and sufficient conditions for when a principal $J$-ideal is a $B(H)$-ideal and we employ these conditions to characterize them (Theorems \ref{T:1}-\ref{T:2}). We also investigate the lattice structure of subideals and principal subideals (building blocks of subideals just as principal ideals are  building blocks for ideals in all rings, see Remark \ref{R:1}(iv)).
Our methods are largely purely algebraic.

Motivated by the advances on ideals of the last decade (for example the semiring structure of the lattice of $B(H)$-ideals, esp. their additive and multiplicative structure, see Remark \ref{R:1}(iv)), we found that bypassing the Lie ideal considerations of Fong-Radjavi \cite[Theorem 1]{FR83},
the positive operator case, we can prove the main theorem \cite[Theorem 2]{FR83} more generally, and we believe more simply and directly.
These advances, to be sure, evolved from works such as \cite{FR83}.
Indeed, the proof of \cite[Lemma 6.3]{DFWW} shares some of the attributes of the proof of \cite[Theorem 1]{FR83},
in particular, the use of a unitary from $H \rightarrow H^{\oplus m}$,
and their use of Proposition \ref{P:3}(i) below for principal $K(H)$-ideals is implicit in their proofs.

\section{PRELIMINARIES}

Every $B(H)$-ideal, $J$, is linear because for each $\alpha \in \mathbb C$, $\alpha I \in B(H)$ so that for each $A \in J$, $(\alpha I)A = \alpha A \in J$.
But a subideal (i.e., a $J$-ideal) may not be linear (Example \ref{E:2}).
This led the authors to introduce linear $J$-ideals in addition to $J$-ideals.
Much the same can be said for real linear $J$-ideals (ideals closed under real scalar multiplication) but we will say little more about these. We start with the following definitions, noting the obvious fact that, intersections of ideals in any ring are themselves ideals.\\

\bD{D:d1}Let $J$ be a $B(H)$-ideal and $S \in J$.

\textbullet \, The principal $B(H)$-ideal generated by the single operator $S$ is given by
\begin{center}
$\left(S\right)$ := $\bigcap \{I \mid I$ is a $B(H)$-ideal containing S$\}$
\end{center}

\textbullet \, The principal linear $J$-ideal generated by $S$ is given by
\begin{center}
$\left(S\right)_J$ := $\bigcap \{\mathcal I \mid \mathcal I$ is a linear $J$-ideal containing $S\}$
\end{center}

\textbullet \, The principal $J$-ideal generated by $S$ is given by
\begin{center}
$\left<S\right>_J$ := $\bigcap \{\mathcal I \mid \mathcal I$ is a $J$-ideal containing $S\}$
\end{center}

\textbullet \, The principal  real linear $J$-ideal generated by $S$ is given by
\begin{center}
$(S)_J^{\mathbb R}$ := $\bigcap \{\mathcal I \mid \mathcal I$ is real linear $J$-ideal containing $S\}$
\end{center}

\noindent To make precise the notion of ideal generation by a set beyond single operator generation, we give the
following natural standard definition.\\

\textbullet \, As above for principal $J$-ideals, likewise for an arbitrary subset $\mathscr S \subseteq J$, $(\mathscr S)$, $(\mathscr S)_J$, $\left<\mathscr S\right>_J$ and $(\mathscr S)_J^{\mathbb R}$ denote respectively, the smallest $B(H)$-ideal, the smallest linear $J$-ideal, the smallest $J$-ideal, and the smallest real linear $J$-ideal generated by the set $\mathscr S$.
\eD

\bR{R:1} \emph{Standard facts on operator ideals.}\\

\nr{i} \cite[Sections 2.8, 4.3]{DFWW} (see also \cite[Section 4]{KW07}): If $I,J$ are $B(H)$-ideals then the product $IJ$, which is both associative and commutative, is the $B(H)$-ideal given by the characteristic set $\Sigma(IJ) = \{\xi \in c_o^* \mid \xi \leq \eta\rho$ for some $\eta \in \Sigma(I)$ and $\rho \in \Sigma(J)\}$.

In abstract rings, the ideal product is defined as the class of finite sums of products of two elements, $IJ := \{\displaystyle{\sum_{finite}} a_ib_i \mid a_i \in I, b_i \in J\}$, but in $B(H)$ the next lemma shows finite sums of operator products defining $IJ$ can be reduced to single products. \\

\nr{ii} \cite[Lemma 6.3]{DFWW} Let $I$ and $J$ be proper ideals of $B(H)$. If $A \in IJ$, then $A = XY$ for some $X \in I$ and $Y \in J$.\\

\nr{iii} \cite[Section 1]{KW07} For $T \in B(H)$, $s(T)$ denotes the sequence of s-numbers of $T$.
Then $ A \in (T)$ if and only if $s(A) = \textrm{O}(D_m(s(T)))$ for some $m \in \mathbb{N}$.
Moreover, $A \in I$ a $B(H)$-ideal if and only if $A^* \in I$ if and only if $|A| \in I$ (via the polar decomposition $A = U|A|$ and $U^*A = |A| = (A^*A)^{1/2}$), with all equivalent to $\diag s(A) \in I$.\\

\nr{iv} The lattice of $B(H)$-ideals forms a commutative semiring with multiplicative identity $B(H)$.
That is, the lattice is commutative and associative under ideal addition and multiplication (see \cite[Section 2.8]{DFWW}) and it is distributive.
Distributivity with multiplier $K(H)$ is stated without proof in \cite[Lemma 5.6-preceding comments]{KW07}.
The general proof is simple and is as follows. For $B(H)$-ideals $I,J,K$, one has
$I,J \subset I+J := \{A+B \mid A \in I, B \in J\}$ and so $IK, JK \subset (I+J)K$, so one has $IK + JK \subset (I+J)K$.
The reverse inclusion follows more simply if one invokes (ii) above: $X \in (I+J)K$ if and only if $X = (A+B)C$ for some $A \in I, B \in J, C \in K$.

The lattice of $B(H)$-ideals is not a ring because, for instance, $\{0\}$ is the only $B(H)$-ideal with an additive inverse, namely, $\{0\}$ itself, so it is not an additive group. It is also clear that $B(H)$ is the multiplicative identity but no $B(H)$-ideal has a multiplicative inverse.

One importance of principal ideals in a general ring is that they are building blocks for all ideals $I$ that contain them in that:
$I \quad = \displaystyle{\bigcup_{r_1,\dots,r_n \in I, \,n \in \mathbb N}} (r_1)+ \cdots +(r_n)$

\nr{v} When $T = \displaystyle{\sum_{i=1}^{n}}A_{i}TB_{i}$ with each $A_i$ or $B_i \in J$, the important s-number relation holds:
$s(T) = \textrm{O}(D_{m}(T)s(C))$ for some $C \in J$
(since then $T \in (T)J$ and so follows from \cite[Section 1, p. 6]{KW07} and Remark \ref{R:1}(i)).\\
\eR

\noindent \emph{Algebraic description of principal subideals of $B(H)$}

\bP{P:3}
(i) For $S \in J$, an algebraic description of principal linear $J$-ideal $(S)_J$, principal real linear $J$-ideal $(S)_J^{\mathbb R}$ and principal $J$-ideal $\left<S\right>_J$ are given by
\begin{center}
$(S)_J = \{\alpha S + AS + SB + \displaystyle{\sum_{i=1}^{m}}A_iSB_i \mid  A,\,B,\,A_i, \,\,B_i \in J,\,\alpha \in \mathbb{C},\,m \in \mathbb{N}\}$
\end{center}
\begin{center}
$(S)_J^{\mathbb R} = \{ rS + AS + SB + \displaystyle{\sum_{i=1}^{m}}A_iSB_i \mid  A,\,B,\,A_i, \,\,B_i \in J,\,r \in \mathbb{R},\,m \in \mathbb{N}\}$
\end{center}
\begin{center}
$\left<S\right>_J = \{ nS + AS + SB + \displaystyle{\sum_{i=1}^{m}}A_iSB_i \mid  A,\,B,\,A_i, \,\,B_i \in J,\, n \in \mathbb{Z},\,m \in \mathbb{N}\}$
\end{center}

\noindent(ii)\quad $J(S)J = \{\displaystyle{\sum_{i=1}^{m}}A_{i}SB_{i} \mid A_{i}, B_{i} \in J$ and $ m \in \mathbb{N}\}$ = $J^{2}(S)$
\eP
\bp \quad \\
(i) Define $\mathfrak{S}:=\{\alpha S + AS + SB + \displaystyle{\sum_{i=1}^{m}}A_iSB_i\mid A, B, A_i, B_i \in J,\alpha\in\mathbb{C},
m \in \mathbb{N}\}$. \\
It is easy to check that $\mathfrak{S}$ is a linear $J$-ideal containing $S$ and that $\mathfrak{S} \subseteq (S)_J$. \\
The reverse inclusion holds since, as an intersection of linear $J$-ideals, $(S)_J$ is the smallest linear $J$-ideal containing $S$.

Similar are the proofs for the forms for $(S)_J^{\mathbb R}$ and $\left<S\right>_J$.\\

\noindent (ii) Clearly $\{\displaystyle{\sum_{i=1}^{m}}A_{i}SB_{i} \mid A_{i}, B_{i} \in J\} \subseteq J(S)J$. For the reverse inclusion,
every element of $J(S)J$ is of the form $AXB$ for $A, B \in J$ and $X \in (S)$ (Remark \ref{R:1}(ii))
and elements of $(S) = (S)_{B(H)}$ are of the form $X = \displaystyle{\sum_{i=1}^{n}}C_iSD_i$ for $C_{i}, D_{i} \in B(H)$.
So $AXB = \displaystyle{\sum_{i=1}^{n}} AC_iSD_iB  \in \{\displaystyle{\sum_{i=1}^{m}}A_{i}SB_{i} \mid A_{i}, B_{i} \in J\}$.
That $J(S)J = J^{2}(S)$ follows from $B(H)$-ideal product commutativity (Remark \ref{R:1}(i)).
\ep

Immediate from this proposition one has the explicit descriptions of $(S)_J $, $\left<S\right>_J$ and $(S)_J^{\mathbb R}$:
\bC{C:1} For $S \in J$,
\begin{equation*}
(S)_J = \mathbb{C}S + JS + SJ + J(S)J
\end{equation*}
\begin{equation*}
\left<S\right>_J = \mathbb{Z}S + JS + SJ + J(S)J
\end{equation*}
\begin{equation*}
(S)_J^{\mathbb R} = \mathbb{R}S + JS + SJ + J(S)J
\end{equation*}
\eC
An extension of an ideal notion called soft-edged is essential for our generalization of Fong-Radjavi's work \cite{FR83}.
Soft-edged ideals (soft ideals for short), that is, $B(H)$-ideals $I$ for which $IK(H) = I$,
were first introduced by Kaftal and Weiss in \cite{KW02} and \cite[Section 3, esp., Definition 3.1]{KW10} and studied further in \cite{KW07} as part of a study on traces
motivated in part by Dixmier's implicit use of softness to construct the so-called Dixmier trace \cite{D66}.
However, as Remark \ref{R:2} below indicates, these softness notions involving $K(H)$ appeared some years earlier.

\bD{D:1} \quad \\
For $B(H)$-ideals $I$ and $J$, the ideal $I$ is called ``$J$-soft" if $IJ = I$. \\
Equivalently in the language of s-numbers (see Remark \ref{R:1}(i),(ii),(v)): \\
For every $A \in I$, $s_n(A) = O(s_n(B)s_n(C))$ for some $B \in I, C \in J, m \in \mathbb N$.
\eD

\noindent Clearly $J$-softness of $I$ implies $I \subseteq J$, so this notion applies only to those $B(H)$-ideals already in $J$.

\begin{remark} \label{R:2}
Recently A. Pietsch alerted us that notions of softness, soft interior and soft complement in particular,
arose years earlier in Banach spaces and called by other names (private communication).
Soft-edged ideals ($K(H)$-soft) are those $B(H)$-ideals $I$ that are equal to their soft interiors $IK(H)$ as defined in \cite[Definition 3.1]{KW10}.
``For Banach ideals over the Hilbert space,
B.S. Mityagin \cite{bM67} introduced the properties of being ``minimal'' and ``maximal''.
These concepts were generalized to arbitrary operator ideals over Banach spaces
by A. Pietsch \cite[4.8.2+6 and 4.9.2+6]{aP80}.
It turned out that because of their local structure,
maximal quasi-Banach operator ideals are of particular importance.
Obviously, the following equations hold:
minimal kernel = soft interior and maximal hull = soft cover.''
(As defined in \cite[Definition 3.1 and succeeding \P]{KW10} and \cite[Definition 4.1]{KW07}, the soft cover of $I$ denoted there $sc\,I$ is the quotient ideal $I/K(H)$.)
\end{remark}
\newpage

\section{PROOFS OF MAIN RESULTS}

We first reduce condition (i) of Theorem \ref{T:1} to the \emph{principal linear $J$-ideal} case via the following lemma.

\bL{L:3.1}
$\left<S\right>_J$ or $(S)_J^{\mathbb{R}}$ is a $B(H)$-ideal if and only if  $(S)_J$ is a $B(H)$-ideal. And in this case, $$ \left<S\right>_J = (S)_J^{\mathbb{R}} = (S)_J = (S).$$
\eL
\bp
$\Rightarrow$: If $\left<S\right>_J$ is a $B(H)$-ideal, then $\frac{1}{2}S \in \left<S\right>_J$.
For principal $J$-ideals, Corollary \ref{C:1} insures that
\begin{equation*}
\frac{1}{2}S = nS + AS + SB + A'SB',
\end{equation*}
hence $(\frac{1}{2}- n)S \in J(S)$.
As $J(S)$ is linear and $\frac{1}{2}\neq n$ one has $S \in J(S)$ in which case, by minimality, $(S) \subset J(S)$,
and since the reverse inclusion is automatic, $(S) = J(S) = J(S)J$
(the second equality follows from the first).
But $J(S)J \subseteq (S)_J$ (Corollary \ref{C:1}) so $(S)\subseteq (S)_J$.
And since $(S)$ is a $B(H)$-ideal, it is also a $J$-ideal, so by minimality, $(S)_J \subseteq (S)$, and hence $(S)=(S)_J$ is a $B(H)$-ideal.

Similarly, for principal real $J$-ideals, if $(S)_J^{\mathbb{R}}$ is a $B(H)$-ideal, then $iS \in (S)_J^{\mathbb{R}}$ where $i = \sqrt{-1}$,
so $iS = rS + AS + SB + A'SB'$, and since $i \neq r \in \mathbb{R}$, likewise $S \in J(S)$, then likewise $(S) = (S)_J$.

That all three are equal is proved next by the reverse implication.

$\Leftarrow$: If $(S)_J$ is a $B(H)$-ideal, then $(S) = J(S)J$ (see proof below for Theorem \ref{T:1} (i) $\Rightarrow$ (ii)).
Then from Corollary \ref{C:1}, $(S)_J \supseteq (S)$ (as above the reverse inclusion is automatic),
and likewise $\left<S\right>_J = (S)$ and $(S)_J^{\mathbb{R}} = (S)$.
Therefore all three are the same $B(H)$-ideal $(S)$.
\ep
From this lemma we can now prove Theorem \ref{T:1} replacing condition (i) with its equivalent that the \textit{principal linear $J$-ideal} $(S)_J$ is a $B(H)$-ideal.\\

\emph{Proof of Theorem \ref{T:1}.}\\

The proof follows the order:
\begin{equation*}
(i)\Rightarrow (ii) \Rightarrow (iii) \Rightarrow (iv) \Rightarrow (ii) \Rightarrow (i).
\end{equation*}
\emph{(i) $\Rightarrow$ (ii)}: For this implication we give two proofs. The first is a primitive proof using s-numbers.
Distilled from this, the second is more modern and shorter reflecting the perspective from advances on ideals from the last decade.
Moreover it is the method (and notation) we use later to generalize this theorem and its preliminary Lemma \ref{L:3.1} to finitely generated $J$-ideals
(Lemma \ref{L:3.2} - Theorem \ref{T:3}).\\

\newpage

\emph{Proof 1.} \\
Denoting $\xi := s(S)$ the sequence of s-numbers of $S$, to show $J(S) = (S)$ it suffices to show $\diag \xi \in  J(S)$
so $S \in J(S)$ (Remark \ref{R:1}(iii)) and then this equality holds as we saw above in the proof of Lemma \ref{L:3.1}.
Since $(S)_{J}$ is a $B(H)$-ideal containing $S$, then $\diag \, (\xi_{1}, 0 , \xi_{3}, 0 , \dots), \diag \, (0 , \xi_{2}, 0 ,\xi_{4}, 0 , \dots) \in (S)_{J}$
by multiplying $\diag s(S)$ by suitable diagonal projections.
Then by Corollary \ref{C:1} for principal linear $J$-ideal $(S)_J$,

(1)\qquad $\diag (\xi_{1}, 0 , \xi_{3}, 0 , \dots) = \alpha S+ AS+ SB+ CS'D$ \\
for some $\alpha \in \mathbb C, A, B, C, D \in J$ and $S' \in (S)$.\\

(2)\qquad  $\diag (0 , \xi_{2}, 0 ,\xi_{4}, 0 , \dots)$ = $\beta S+ A'S+ SB'+ C'S''D'$\\
for some $\beta \in \mathbb C, A', B', C', D' \in J$ and $S'' \in (S)$.\\

\noindent If $\alpha$ = 0 or $\beta = 0$, then $\diag \, (\xi_{1}, 0 , \xi_{3}, 0 , \dots)$ or $\diag \, (0 , \xi_{2}, 0 ,\xi_{4}, 0 , \dots) \in J(S)$.
Then, in either case, because  $(\xi_{3}, \xi_{5}, \dots) \le (\xi_{2}, \xi_{4}, \dots) \le (\xi_{1}, \xi_{3}, \dots)$,
one has $\diag \xi \in J(S)$ by the hereditariness of $\Sigma(J(S))$,
hence $S \in J(S)$ which again, as in the proof of Lemma \ref{L:3.1}, is equivalent to condition (ii): $(S) = J(S)$.
So without loss of generality assume $\alpha, \beta \ne 0$. Multiplying (1) by $-\beta$, (2) by $\alpha$ and adding,
$\diag \, (-\beta\xi_{1}, \alpha\xi_{2} , -\beta\xi_{3}, \alpha\xi_{4}, \dots) = (-\beta A + \alpha A')S + S(-\beta B + \alpha B') + CS'D + C'S''D' \in J(S).$
\noindent Again multiplying by suitable diagonal projections,\\
$\diag \, (|\beta| \xi_{1}, 0 , |\beta|\xi_{3}, 0, \dots) \in J(S)$, and so likewise $\diag \, (\xi_{1}, 0, \xi_{3}, 0, \dots)$ and then also $\diag \, (\xi_{1}, \xi_{3}, \dots)$.
By hereditariness $\diag \, (0 , \xi_{2}, 0 ,\xi_{4}, 0 , \dots)\in J(S)$ because $(\xi_{2}, \xi_{4}, \dots) \le (\xi_{1}, \xi_{3}, \dots)$.
Hence $\diag \xi \in J(S)$ which again is equivalent to condition (ii).\\

\textit{Proof 2.} \\
For any unitary map $\phi: H \rightarrow H \oplus H$, $S \rightarrow \phi S \phi^{-1}$ preserves s-number sequences and ideals via Calkin's representation.
Since $(S)_{J}$ is a $B(H)$-ideal containing $S$, $\phi^{-1}(S \oplus 0)\phi,\, \phi^{-1}(0 \oplus S)\phi \in (S)_{J}$ since they possess the same s-numbers as $S$.
Then by Corollary \ref{C:1} for principal linear $J$-ideal $(S)_J$,\\

(3)\qquad $\phi^{-1}(S \oplus 0)\phi = \alpha S+ AS+ SB+ CS'D$\\
for some $\alpha \in \mathbb C,~A, B, C, D \in J$ and $S' \in (S)$.\\

(4)\qquad $\phi^{-1}(0 \oplus S)\phi = \beta S+ A'S+ SB'+ C'S''D'$ \\
for some $\beta \in \mathbb C,~A', B', C', D' \in J$ and $S'' \in (S)$.\\

\noindent If $\alpha$ = 0 or $\beta = 0$, then $\phi^{-1}(S \oplus 0)\phi$ or $\phi^{-1}(0 \oplus S)\phi \in J(S)$.
Then, in either case, $S \in J(S)$ which, as we saw in the proof of Lemma \ref{L:3.1}, is equivalent to condition (ii): $(S) = J(S)$.
Finally if $\alpha, \beta \ne 0$, multiplying (3) by $-\beta$, (4) by $\alpha$ and adding obtains
$\phi^{-1}(-\beta S \oplus \alpha S)\phi = (-\beta A + \alpha A')S + S(-\beta B + \alpha B') + CS'D + C'S''D'$ which belong to $J(S)$.
Multiplying in $B(H \oplus H)$ by a suitable diagonal projection one obtains $\phi^{-1}(S \oplus 0)\phi \in J(S)$.
Hence, also $S \in J(S)$, again equivalent to (ii).

\emph{(ii) $\Rightarrow$ (iii)}:

\noindent If $(S)$ = $J(S)$, since $B(H)$-ideals commute and their semiring multiplication is associative,\\
$J(S)J = \{J(S)\}J = (S)J = J(S) = (S)$.
Therefore, $S= X\left(\displaystyle{\sum_{i=1}^{k}}C_{i}SD_{i}\right)Y$, for some $X, Y \in J$ and $C_{i}, D_{i} \in B(H)$ (\cite[Lemma 6.3]{DFWW}), hence (iii):\\
$S= \displaystyle{\sum_{i=1}^{k}}(XC_{i})S(D_{i}Y) \in J(S)J$.\\

\emph{(iii) $\Rightarrow$ (iv)}: This is Remark \ref{R:1}(v).\\

\emph{(iv) $\Rightarrow$ (ii)}:\\
This follows directly from the definition of the characteristic set of the product $\Sigma ((S)J)$ and then that $s(S) = \textrm{O}(D_{k}(s(S))s(T))$ for some $T \in J$ implies $\diag \, (s(S)) \in (S)J$. Therefore, $S \in (S)J$ and (ii) follows.\\

\emph{(ii) $\Rightarrow$ (i)}:\\
$(S) = (S)J$ implies $(S) = J(S)J$. Therefore, $S = \displaystyle{\sum_{i=1}^{n}}A_{i}SB_{i}$ for some $A_{i}$ and $B_{i}$ in $J$.
By Corollary \ref{C:1} $(S)_J = \mathbb{C}S + JS + SJ + J(S)J$ so $(S)_J \subseteq J(S)J$, and substituting $S$ by $\displaystyle{\sum_{i=1}^{n}}A_{i}SB_{i}$ in the right-hand side equation  one obtains $J(S)J \subseteq (S)_J$,\\
so $(S)_J = J(S)J$.
Since $J(S)J$ is a $B(H)$-ideal, $(S)_J$ is a $B(H)$-ideal.
Moreover, $(S)_J = (S) = J(S)J$.\\

\noindent Remark on Theorem \ref{T:1}--Proof (i)$\Rightarrow$(ii):
Proof 2 may seem simpler or shorter but Proof 1 keeps the analysis in the same Hilbert space, it is more constructive, and it appears to us more useful.\\

\emph{Proof of Theorem \ref{T:2}}\\

Corollary \ref{C:1} gives directly the explicit descriptions of $(S)_J $, $\left<S\right>_J$ and $(S)_J^{\mathbb R}$:
\begin{equation*}
\left<S\right>_J = \mathbb Z S + JS + SJ + J(S)J
\end{equation*}
\begin{equation*}
(S)_J = \mathbb C S + JS + SJ + J(S)J
\end{equation*}
\begin{equation*}
(S)_J^\mathbb R = \mathbb R S + JS + SJ + J(S)J
\end{equation*}
from which it follows that
\begin{center}
\begin{equation*}
J(S)J \subseteq JS + SJ + J(S)J \subseteq \left<S\right>_J \subseteq (S)_J^{\mathbb R} \subseteq (S)_J \subseteq (S)
\end{equation*}
\end{center}
That $J(S)J$ is a $B(H)$-ideal and $JS + SJ + J(S)J$ is a $J$-ideal is clear.
An immediate consequence of Lemma \ref{L:3.1} is that $J$-softness of $(S)$ is equivalent to having at least one and hence all three principal $J$-ideals be the $B(H)$-ideal $(S)$.

\bC{C:2}
Equality of at least one of the inclusions in $$\left<S\right>_J \subseteq (S)_J^{\mathbb R} \subseteq (S)_J \subseteq (S)$$ implies equality throughout.
In that case, all three types of principal $J$-ideals are the $B(H)$-ideal $(S)$.
\eC
\bp
If $\left<S\right>_J = (S)_J^{\mathbb R}$, then $\frac{1}{2}S \in \left<S\right>_J$, and if $(S)_J^{\mathbb R} = (S)_J$, then $iS \in (S)_J^{\mathbb R},$\\
$i = \sqrt{-1}$.
By Theorem \ref{T:2}, $\frac{1}{2}S = nS + AS + SB + A'SB'$ for some $n \in \mathbb N$ in first case and $iS = rS + AS + SB + A'SB'$ for some $r \in \mathbb R$ in the second case.
Since $\frac{1}{2} \neq n$ and $i \ne r$, one has $S \in J(S)$.
So for either set equality, $(S)$ is $J$-soft which, again by Theorem \ref{T:2}, implies $J(S) = (S) = \left<S\right>_J = (S)_J = (S)_J^\mathbb R$.
For the last equality, then $(S)_J$ is the $B(H)$-ideal $(S)$ so Lemma \ref{L:3.1} applies.
\ep
Corollary \ref{C:2} provides easy examples where the inclusions of the three types of principal $J$-ideals are proper, that is
all principal $B(H)$-ideals that are not $J$-soft as in the following example.

\bE{E:1}
For $J = K(H)$ and $S = \diag\left<\frac{1}{n}\right>$, $(\diag\left<\frac{1}{n}\right>)$ is not
$K(H)$-soft.
If it were, then $(\diag\left<\frac{1}{n}\right>) = (\diag\left<\frac{1}{n}\right>)K(H)$
which further implies\\ $\left<\frac{1}{n}\right> \in \Sigma ((\diag\left<\frac{1}{n}\right>)K(H))$.
By Remark \ref{R:1}(i),\,(iii) and $\Sigma(K(H)) = c_{0}^{*}$, one has
$\left<\frac{1}{n}\right> = o(D_{m}\left<\frac{1}{n}\right>)$ for some $m \in \mathbb N$, contradicting
$\left(\frac{\left<\frac{1}{n}\right>}{D_{m}\left<\frac{1}{n}\right>}\right)_k = \frac{\frac{1}{mj+r}}{\frac{1}{j+1}} \rightarrow \frac{1}{m}$
as $k \rightarrow \infty$ where $k = mj+r$.
But since $(\diag\left<\frac{1}{n}\right>)$ is not $K(H)$-soft, by Theorem \ref{T:2}, the three inclusions above are strict.
\eE

Finishing the discussion on inclusions, $J(S)J \subseteq JS + SJ + J(S)J$ can be proper as given below in Example \ref{E:2}.
\bE{E:2}
If $S = \diag \left<\frac{1}{n}\right>$ and $J = (S)$,
then $J(S)J \subsetneq J(S)J + JS + SJ$ because $\diag \left<\frac{1}{n^2}\right> \in JS \setminus J(S)J$.
Assume otherwise that $\diag \left<\frac{1}{n^2}\right> \in J(S)J.$
It is elementary to show, via characteristic sets using Remark \ref{R:1}(iii),\\
that $(A)(B) = (AB)$ when $A$ and $B$ are simultaneously diagonalizable with entries of the diagonalized operators in non-increasing order.\\ Therefore, $J(S)J =(S)(S)(S)=(S^{3})$.
Then again by Remark \ref{R:1}(iii) one obtains $\left<\frac{1}{n^2}\right> = \textrm{O}(D_{m}\left<\frac{1}{n^{3}}\right>)$ thereby contradicting $\left(\frac{\left<\frac{1}{n^{2}}\right>}{D_{m}\left<\frac{1}{n^{3}}\right>}\right)_{k} = \frac{\frac{1}{(mj+r)^{2}}}{\frac{1}{(j+1)^{3}}} \rightarrow \infty$ as $k \rightarrow \infty$, where $k = mj+r$.
So the inclusion is strict.
\eE

\newpage

\bE{E:3}\emph{A concrete nonlinear principal ideal: $\left<\diag\left<\frac{1}{n}\right>\right>_{K(H)}$}\\
Example \ref{E:1} also provides us with a concrete principal nonlinear $K(H)$-ideal and leads to a plethora of them.
The principal $K(H)$-ideal, $\left<\diag\left<\frac{1}{n}\right>\right>_{K(H)}$ is not linear. If it were linear,
then $i\diag\left<\frac{1}{n}\right> \in \left<\diag\left<\frac{1}{n}\right>\right>_{K(H)}$.
By Corollary \ref{C:1}, $i\diag\left<\frac{1}{n}\right> = m\diag\left<\frac{1}{n}\right> + A \diag\left<\frac{1}{n}\right> + \diag\left<\frac{1}{n}\right>B + \displaystyle{\sum_{i=1}^{n}A_{i}\diag\left<\frac{1}{n}\right>B_{i}}$ for some $A, B, A_{i}, B_{i} \in K(H), m \in \mathbb{N}$.
So $(i-m)\diag\left<\frac{1}{n}\right> \in (\diag\left<\frac{1}{n}\right>)K(H)$ and since $i \neq m$,
one has $\diag\left<\frac{1}{n}\right> \in (\diag\left<\frac{1}{n}\right>)K(H)$ hence $(\diag\left<\frac{1}{n}\right>)$ is $K(H)$-soft contradicting Example \ref{E:1}.
\eE

\vspace{-.15in}

\section{FINITELY GENERATED $J$-IDEALS}

Theorems \ref{T:1}-\ref{T:2} and its preliminary Lemma \ref{L:3.1} generalize to finitely generated $J$-ideals somewhat naturally in Lemma \ref{L:3.2} and Theorems \ref{T:3}-\ref{T:4} below.

Moreover, every finitely generated $B(H)$-ideal is always a principal $B(H)$-ideal because, as is straightforward to see,
the $B(H)$-ideal generated by\\ $\mathscr S = \{S_1, \dots, S_n\} \subset B(H)$ is the principal ideal $(\mathscr S) = (|S_1| + \dots + |S_n|)$ (recall Definition \ref{D:d1} and that $|S| := (S^*S)^{1/2}$).
But finitely generated $J$-ideals (general, linear or real linear) may not be principal as seen in the following example.

\bE{E:4} \emph{[A nonprincipal doubly generated $J$-ideal of any of the three types]}

For $J=K(H)$, $S_1 = \diag \, (1, 0, \frac{1}{2}, 0, \frac{1}{3}, \cdots)$ and $S_2 = \diag \, (0, 1, 0, \frac{1}{2}, 0, \frac{1}{3}, \cdots)$,
$(S_1,S_2)_{K(H)}$ is not a principal linear $K(H)$-ideal, and likewise for the general and real linear cases $\left<S_{1},S_{2}\right>_J$ and $(S_{1},S_{2})_J^{\mathbb R}$.

We omit our proof of Example \ref{E:4} because it led to a generalization providing, under certain assumptions, a necessary and sufficient condition for a $J$-ideal generated by two operators to be a principal $J$-ideal (see Proposition \ref{P:7.2}).
\eE

To investigate finitely generated $J$-ideals we start with algebraic descriptions of the three types of $J$-ideals with $N$ generators analogous to Proposition \ref{P:3}(i).
Observing that $(\{S_{1}, \cdots, S_{N}\})_J = (S_{1})_J + \cdots + (S_{N})_J$ and similarly for the other two types (this holds for ideals in general rings),
and using the same arguments used for the algebraic description of principal $J$-ideals in Proposition \ref{P:3}(i), it is straightforward to see that denoting
$\mathscr S := \{S_{1}, \cdots, S_{N}\}$,

\bP{P:4} For $\mathscr S \subseteq J$, $J(\mathscr S)J = J(|S_{1}| + \cdots + |S_{N}|)J$ and
\begin{equation*}
(\mathscr S)_J = \mathbb{C}S_{1} + \cdots + \mathbb{C}S_{N} + JS_{1} + \cdots + JS_{N} +  S_{1}J + \cdots + S_{N}J + J(\mathscr S)J
\end{equation*}
\begin{equation*}
\left<\mathscr S\right>_J = \mathbb{Z}S_{1} + \cdots + \mathbb{Z}S_{N} + JS_{1} + \cdots + JS_{N} +  S_{1}J + \cdots + S_{N}J + J(\mathscr S)J
\end{equation*}
\begin{equation*}
(\mathscr S)_J^{\mathbb R} = \mathbb{R}S_{1} + \cdots + \mathbb R S_{N} + JS_{1} + \cdots + JS_{N} +  S_{1}J + \cdots + S_{N}J + J(\mathscr S)J
\end{equation*}
\eP

As we first reduced condition (i) of Theorem \ref{T:1} to the \textit{linear $J$-ideal} case,
now we reduce condition (i) of the next Theorem \ref{T:3} to the \textit{linear $J$-ideal} case via the following lemma.

\bL{L:3.2} For $\mathscr S := \{S_{1}, \cdots, S_{N}\} \subseteq J$,
$\left<\mathscr S\right>_J$ or $(\mathscr S)_J^{\mathbb{R}}$ is a $B(H)$-ideal if and only if  $(\mathscr S)_J$ is a $B(H)$-ideal.
In this case, they are all the $B(H)$-ideal spanned by $\mathscr S$, that is,
\begin{equation*}
\left<\mathscr S\right>_J = (\mathscr S)_J^{\mathbb{R}} = (\mathscr S)_J = (\mathscr S) = (|S_{1}|+ \cdots + |S_{N}|).
\end{equation*}
\eL
\bp
($\Rightarrow$:)\\
For any unitary $\phi: H \rightarrow \underbrace{H \oplus \cdots \oplus  H}_{N+1-times}$, recall as mentioned earlier that the map $S \rightarrow \phi S\phi^{-1}$ preserves s-numbers and ideals (via Calkin's representation).
That $\left<\mathscr S\right>_J$ is a $B(H)$-ideal implies that the $N+1$ operators in $B(H)$,\\
$\phi^{-1}(S_{1} \oplus \underbrace{0 \oplus \cdots \oplus 0}_{N-times})\phi, \cdots, \phi^{-1}(\underbrace{0 \oplus \cdots \oplus 0}_{N-times} \oplus S_{1})\phi \in \left<S_{1}, \cdots, S_{N}\right>_J$.
By Proposition \ref{P:4} one obtains the set of $N+1$ equations:
\begin{align*}
\phi^{-1}(S_{1} \oplus \underbrace{0 \oplus \cdots \oplus 0}_{N-times})\phi &= n_{1,1}S_{1} + \cdots + n_{1,N}S_{N} + A_{1}\\
\vdots\\
\phi^{-1}(\underbrace{0 \oplus \cdots \oplus 0}_{N-times} \oplus S_{1})\phi &= n_{N+1,1}S_{1} + \cdots + n_{N+1,N}S_{N} + A_{N+1}
\end{align*}
where $A_{j} \in J(\mathscr S)$ for $1 \leq j \leq N+1$ and $n_{i,j} \in \mathbb Z$ for $1 \leq i \leq N+1, 1 \leq j \leq N$. \\
By basic linear algebra, the $N+1$ vectors in $\mathbb C^{N}$:\\ $(n_{1,1}, \cdots, n_{1,N}),\,  \cdots, (n_{N+1,1}, \cdots, n_{N+1,N})$ are  linearly dependent.
That is, for some $\alpha_{1}, \cdots, \alpha_{N+1}$ not all $0$, say $\alpha_k \ne 0$, and one has
\begin{equation*}
\alpha_{1}\left<n_{1,1}, \cdots, n_{1,N}\right> + \cdots + \alpha_{N+1}\left<n_{N+1,1}, \cdots, n_{N+1,N}\right> = \left<0, \cdots, 0\right>
\end{equation*}
and hence,\\
$\alpha_{1}\phi^{-1}(S_{1} \oplus \underbrace{0 \oplus \cdots \oplus 0}_{N-times})\phi + \cdots + \alpha_{N+1}\phi^{-1}(\underbrace{0 \oplus \cdots \oplus 0}_{N-times} \oplus S_{1})\phi\\
= \alpha_{1} A_1 + \cdots + \alpha_{N+1}A_{N+1}$ which is in $J(\mathscr S)$, that is,
\begin{equation*}
\phi^{-1}(\alpha_{1}S_{1} \oplus \alpha_{2}S_{1} \oplus \cdots \oplus \alpha_{N+1}S_{1})\phi \in J(\mathscr S).
\end{equation*}
Multiplying $\alpha_{1}S_{1} \oplus \alpha_{2}S_{1} + \cdots + \alpha_{N+1}S_{1}$ by a suitable diagonal projection one obtains
\begin{equation*}
\phi^{-1}(0 \oplus 0 \cdots \alpha_{k}S_{1} \oplus \cdots \oplus 0)\phi  \in J(\mathscr S)
\end{equation*}
and because $s(\phi^{-1}(0 \oplus 0 \cdots \alpha_{k}S_{1} \oplus \cdots \oplus 0)\phi) = |\alpha_k|s(S_1)$, one obtains $S_1 \in J(\mathscr S)$.
Likewise all $S_j \in J(\mathscr S)$ ($1 \le j \le N$) and hence $(\mathscr S) \subset J(\mathscr S)$.
The reverse inclusion is automatic so $(\mathscr S)$ is $J$-soft, that is, $(\mathscr S) = J(\mathscr S)$.\\
But then $(\mathscr S) = J(\mathscr S) = (\mathscr S)J = J(\mathscr S)J \subseteq (\mathscr S)_J$, and since by minimality $(\mathscr S)_J \subseteq (\mathscr S)$,
one obtains $(\mathscr S)_J = (\mathscr S)$ hence it is a $B(H)$-ideal.

If instead of assuming $\left<\mathscr S \right>_J$ is a $B(H)$-ideal, one assumes that $(\mathscr S)_J^{\mathbb R}$ is a $B(H)$-ideal,
then the proof is essentially the same except for the system of equations where instead of choosing integer scalars $n_{ij} \in \mathbb Z$
one chooses real scalars $r_{ij} \in \mathbb R$.

That all three are equal is proved next by the reverse implication.

$\Leftarrow$: If $(\mathscr S)_J$ is a $B(H)$-ideal, then $(\mathscr S) = J(\mathscr S)J$ (see proof below for Theorem \ref{T:3} (i) $\Rightarrow$ (ii)).
Then from Proposition \ref{P:4}, $(\mathscr S)_J \supseteq (\mathscr S)$ hence $(\mathscr S)_J = (\mathscr S)$ (as again by minimality the reverse inclusion is automatic),
and likewise $\left<\mathscr S\right>_J = (\mathscr S)$ and $(\mathscr S)_J^{\mathbb{R}} = (\mathscr S)$.
And since $(\mathscr S) = (|S_{1}| + \cdots + |S_{N}|)$, all three are equal to the principal $B(H)$-ideal $(|S_{1}| + \cdots + |S_{N}|)$.
\ep

\bR{R:3}
Lemma \ref{L:3.2} implies that if any of the three types of finitely generated $J$-ideals is a $B(H)$-ideal then all three are equal to the principal $B(H)$-ideal generated by the generators.
\eR

The next theorem is the finitely generated analog of Theorem \ref{T:1}.

\bT{T:3}
For $\mathscr S := \{S_{1}, \cdots, S_{N}\} \subseteq J$, the following are equivalent.\\

\nr{i} Any of the 3 types of $J$-ideals generated by $\mathscr S$; $(\mathscr S)_J$, $\left<\mathscr S\right>_J$ or $(\mathscr S)_J^{\mathbb{R}}$, is a $B(H)$-ideal.\\

\nr{ii} The $B(H)$-ideal $(\mathscr S)$ is $J$-soft, i.e., $(\mathscr S)$ = $J(\mathscr S)$ (equivalently, $(\mathscr S)$ = $(\mathscr S)J$).\\

\nr{iii} For all \,$1 \leq j \leq N$,
\begin{equation*}
S_{j} = \displaystyle{\sum_{i=1}^{n_{j}}}A_{i1}S_{1}B_{i1} + \cdots + \displaystyle{\sum_{i=1}^{k_{j}}}A_{iN}^{'}S_{N}B_{iN}^{'}
\end{equation*}
\qquad \quad \text{for some} $A_{ij}, B_{ij}, A_{ij}^{'}, B_{ij}^{'} \in J, n_{j}, k_{j} \in \mathbb{N}$.\\

\nr{iv}  For all \,$1 \leq j \leq N$,
\begin{equation*}
s(S_{j}) = \textrm{O}(D_{m}(s(|S_{1}|+ \cdots +|S_{N}|))s(T)) \text{~for some~} T \in J \text{~and~} m \in \mathbb{N}.\\
\end{equation*}
\eT

\bp
In view of Lemma \ref{L:3.2}, to prove Theorem \ref{T:3} one can replace condition (i) with its equivalent:
$(\mathscr S)_J$ is a $B(H)$-ideal.

The proof follows the order:
\begin{equation*}
(i)\Rightarrow (ii) \Rightarrow (iii) \Rightarrow (iv) \Rightarrow (ii) \Rightarrow (i).
\end{equation*}

\emph{(i) $\Rightarrow$ (ii)}: As just stated, it suffices to assume that the linear $J$-ideal $(\mathscr S)_J$ is a $B(H)$-ideal,
and for this it suffices to prove $(\mathscr S)_J = (\mathscr S)$ and as before, merely to show that $S_{1}, \cdots, S_{N} \in (\mathscr S)J$.
Use the proof of Lemma \ref{L:3.2} ((i) $\Rightarrow$ (ii)) except where the scalars $n_{ij} \in \mathbb Z$ appear in the system of equations, use instead complex coefficients $c_{ij} \in \mathbb C$.\\

\emph{(ii) $\Rightarrow$ (iii)}: If $(\mathscr S)$ = $J(\mathscr S)$, since $B(H)$-ideals commute (Remark \ref{R:1}(i)), then $S_{1}, \cdots, S_{N} \in J(\mathscr S)J$ which by Remark \ref{R:1}(ii) further gives
\begin{align*}
S_{j} &= A(\displaystyle{\sum_{i=1}^{n_{j}}C_{i,j}S_{1}D_{i,j}}+\cdots + \displaystyle{\sum_{i=1}^{k_{j}}C_{i,j}^{'}S_{N}D_{i,j}^{'}})B \\
&= \displaystyle{\sum_{i=1}^{n_{j}}}A_{i,j}S_{1}B_{i,j} + \cdots + \displaystyle{\sum_{i=1}^{k_{j}}}A_{i,j}^{'}S_{N}B_{i,j}^{'}
\end{align*}
\text{for some} $A, B, A_{i,j}, B_{i,j}, A_{i,j}^{'}, B_{i,j}^{'} \in J, n_{j}, k_{j} \in \mathbb{N}~\text{for all}~ 1 \leq j \leq N$.\\

\emph{(iii) $\Rightarrow$ (iv)}: \emph{(iii)} implies $S_{1}, \cdots, S_{N} \in J(\mathscr S)$. Since $(\mathscr S) = (|S_{1}| + \cdots + |S_{N}|)$ and using Remark \ref{R:1}(iii), one gets directly for all $1 \leq j \leq N$: \\
$s(S_{j}) = \textrm{O}(D_{m_{j}}(s(|S_{1}|+ \cdots +|S_{N}|))s(T_{j})) \,\text{for some}\, T_{j} \in J~~\text{and}~ m_{j} \in \mathbb{N}.$
Choosing $m = max \{m_j\}$ and $T = |T_1| + \cdots + |T_N|$ suffices to obtain (iv).\\

\emph{(iv) $\Rightarrow$ (ii)}: It follows directly from the definition of the characteristic set $\Sigma(J(\mathscr S))$ that, for all $1 \leq j \leq N$,
\begin{align*}
S_{j} &= \displaystyle{\sum_{i=1}^{n_{j}}}A_{i,j}S_{1} + \cdots + \displaystyle{\sum_{i=1}^{k_{j}}}A_{i,j}^{'}S_{N} \in J(\mathscr S),
\text{for some}~ A_{i,j}, A_{i,j}^{'},  \in J, n_{j}, k_{j} \in \mathbb{N},
\end{align*}
that is, $S_{j} \in (\mathscr S)J$, so again from which (ii) follows.\\

\emph{(ii) $\Rightarrow$ (i):}
$(\mathscr S) = (\mathscr S)J$ implies $(\mathscr S) = J(\mathscr S)J$.
Therefore, for all $1 \leq j \leq N$,\\
\indent \qquad \qquad \qquad $S_{j} = \displaystyle{\sum_{i=1}^{n_{j}}}A_{i,j}S_{1}B_{i,j} + \cdots + \displaystyle{\sum_{i=1}^{k_{j}}}A_{i,j}^{'}S_{N}B_{i,j}^{'}\hspace{2cm}(*)$\\
\text{for some}~ $A_{i,j}, B_{i,j}, A_{i,j}^{'}, B_{i,j}^{'} \in J, n_{j}, k_{j} \in \mathbb{N}$.
Then $(\mathscr S)_J$ is a $B(H)$-ideal because substituting all $S_{j}$ by $(*)$ on the right side of the equality in
\begin{equation*}
(\mathscr S)_J = \mathbb{C}S_{1} + \cdots + \mathbb{C}S_{N} + JS_{1} + \cdots + JS_{N} +  S_{1}J + \cdots + S_{N}J + J(\mathscr S)J,
\end{equation*}
one obtains $(\mathscr S)_J \subseteq J(\mathscr S)J$ with $(\mathscr S)_J \supseteq J(\mathscr S)J$ automatic (Proposition \ref{P:4}).\ep

Theorem \ref{T:3} combined with Proposition \ref{P:4}, Lemma \ref{L:3.2} and Remark \ref{R:3} provide naturally a finitely generated partial analog for Theorem \ref{T:2}.

\bT{T:4}
In addition to the forms of the three types of finitely generated $J$-ideals generated by $\mathscr S$ given by Proposition \ref{P:4} and to the equivalences given by Lemma \ref{L:3.2} on when any of them is a $B(H)$-ideal, one has
\begin{equation*}
J(\mathscr S)J \subseteq JS_{1} + \cdots + JS_{N} +  S_{1}J + \cdots + S_{N}J + J(\mathscr S)J \subseteq \left<\mathscr S\right>_J \subseteq (\mathscr S)_J^\mathbb R \subseteq (\mathscr S)_J \subseteq (\mathscr S)
\end{equation*}
which first two, $J(\mathscr S)J$ and $JS_{1} + \cdots + JS_{N} +  S_{1}J + \cdots + S_{N}J + J(\mathscr S)J$ respectively, are a $B(H)$-ideal and a $J$-ideal.
Consequently, each of these three kinds of finitely generated $J$-ideals have the common $B(H)$-ideal ``nucleus'' $J(\mathscr S)J$,
with the common $J$-ideal \\
$JS_{1} + \cdots + JS_{N} +  S_{1}J + \cdots + S_{N}J + J(\mathscr S)J$ containing it.

All these finitely generated $J$-ideals, $\left<\mathscr S\right>_J \subseteq (\mathscr S)_J^\mathbb R  \subseteq (\mathscr S)_J$,
collapse to merely
\begin{equation*}
J(\mathscr S)J = (\mathscr S) = \left<\mathscr S\right>_J = (\mathscr S)_J = (\mathscr S)_J^\mathbb R
\end{equation*}
if the finitely generated $B(H)$-ideal $(\mathscr S)$ is $J$-soft, that is, $(\mathscr S) = J(\mathscr S)$,
or if any, and hence all of them, is a $B(H)$-ideal.
Moreover, if $\left<\mathscr S\right>_J = (\mathscr S)_J^\mathbb R$ or $\left<\mathscr S\right>_J = (\mathscr S)_J$, then $(\mathscr S)$ is $J$-soft.
\eT
\bp
In view of Theorem \ref{T:3} combined with Proposition \ref{P:4}, Lemma \ref{L:3.2} and Remark \ref{R:3}, the only thing left to prove is:
if $\left<\mathscr S\right>_J = (\mathscr S)_J^\mathbb R$ or $\left<\mathscr S\right>_J = (\mathscr S)_J$
(which itself implies $\left<\mathscr S\right>_J = (\mathscr S)_J^\mathbb R$), then $(\mathscr S)$ is $J$-soft.
Suppose then that without loss of generality $\left<\mathscr S\right>_J = (\mathscr S)_J^{\mathbb R}$.
Denote $(\mathscr S)_J^0 := JS_{1} + \cdots + JS_{N} +  S_{1}J + \cdots + S_{N}J + J(\mathscr S)J$.
Then by Proposition \ref{P:4} one obtains the following set of $N$-linear equations:
\begin{equation*}
\frac{1}{2}S_{1} = n_{1,1}S_{1} + n_{1,2}S_{2} +  \cdots + n_{1,N}S_{N} + X_{1}
\end{equation*}
\begin{equation*}
\vdots
\end{equation*}
\begin{equation*}
\frac{1}{2}S_{N} = n_{N,1}S_{1} + \cdots + n_{N,N-1}S_{N-1} + n_{N,N}S_{N} + X_{N},
\end{equation*}
where $X_1, \dots X_N \in (\mathscr S)_J^0$,
or equivalently,
\begin{equation*}
(1 - 2n_{1,1})S_{1} - \cdots - 2n_{1,N}S_{N} = 2X_{1}
\end{equation*}
\begin{equation*}
\vdots
\end{equation*}
\begin{equation*}
-2n_{N,1}S_{1} - \cdots + (1 - 2n_{N,N})S_{N} = 2X_{N}.
\end{equation*}
In the quotient space $(\mathscr S)_J^{\mathbb R}/(\mathscr S)_J^0$, these equations become a set of $N$-linear equations given by
\begin{equation*}
(1 - 2n_{1,1})S_{1} - \cdots - 2n_{1,N}S_{N} = 0
\end{equation*}
\begin{equation*}
\vdots
\end{equation*}
\begin{equation*}
-2n_{N,1}S_{1} - \cdots + (1 - 2n_{N,N})S_{N} = 0,
\end{equation*}
where $S_{1}, \cdots, S_{N}$ are now considered as cosets in $(\mathscr S)_J^{\mathbb R}/(\mathscr S)_J^0$.
So one obtains a set of $N$-linear equations in $N$-variables namely $S_{1}, \cdots, S_{N}$ where the integer coefficient matrix is given by
\begin{equation*}
\begin{pmatrix}
1-2n_{1,1}  &  -2n_{1,2}            & \dots           & -2n_{1,N}  \\
-2n_{2,1}             &  1-2n_{2,2} & \dots       & -2n_{2,N}   \\
\vdots               &  \vdots              & \ddots          & \vdots \\
-2n_{N,1}             &  -2n_{N,2}            & \dots           & 1-2n_{N,N} \\
\end{pmatrix}
\end{equation*}
The determinant of this matrix is odd, hence nonzero, because it is the sum of even numbers and one odd number (the product of the diagonal odd entries).
Therefore the matrix is invertible implying that this system of $N$-linear equations has only the trivial solution, i.e., the cosets $S_{j} = 0$ for all $1\leq j \leq N$. Hence the operators $S_{j} \in (\mathscr S)_J^0 \subset J(\mathscr S)$ for all $1\leq j \leq N$ and therefore $(\mathscr S) = J(\mathscr S)$.
\ep
Whether or not $(\mathscr S)_J^\mathbb R = (\mathscr S)_J$ implies $(\mathscr S)$ is $J$-soft remains an open question (Section \ref{Questions}, Question 7).

\section{COMPARISON OF SUBIDEALS WITH $B(H)$-IDEALS : ADJOINTS AND ABSOLUTE VALUES}\label{S: 4}

For $T \in B(H), (T) = (|T|)$, but this need not be true for principal linear $K(H)$-ideals (Example \ref{E:(T) ne (|T|)}).
Moreover, all $B(H)$-ideals are selfadjoint, but not necessarily for principal linear $K(H)$-ideals (Example \ref{E: J ne J^*}).

\bE{E:(T) ne (|T|)}

If $J = K(H)$ and operator $T = \diag \left<\frac{(i)^n}{n}\right>$ so $|T| = \diag \left<\frac{1}{n}\right>$, then $(T)_{K(H)} \neq (|T|)_{K(H)}$.\\
In fact, neither $(|T|)_{K(H)}$ $\subsetneq$ $(T)_{K(H)}$ nor $(T)_{K(H)}$ $\subsetneq$ $(|T|)_{K(H)}$. Indeed, suppose $T \in (|T|)_{K(H)}$ = $(\diag\left<\frac{1}{n}\right>)_{K(H)}$.
By Proposition \ref{P:3},
\begin{equation*}
T = \alpha\, \diag \left<\frac{1}{n}\right> + A \diag \left<\frac{1}{n}\right> + \diag \left<\frac{1}{n}\right>B + \sum_{i=1}^{n}A_{i} \diag\left<\frac{1}{n}\right>B_{i}
\end{equation*}
\text{for some} $A, B, A_{i}, B_{i} \in K(H), n \in \mathbb N.$
\noindent Therefore, $\diag \left<\frac{(i)^{n}-\alpha}{n}\right> \in (\diag \left<\frac{1}{n}\right>)K(H)$ implying
$\diag \left<\frac{|(i)^{n}-\alpha|}{n}\right> \in (\diag \left<\frac{1}{n}\right>)K(H).$\\
Hence, $\diag \left<\frac{1}{2n-1}\right> \in (\diag \left<\frac{1}{n}\right>)K(H)$ implying $\diag \left<\frac{1}{n}\right> \in (\diag \left<\frac{1}{n}\right>)K(H)$,
contradicting Example \ref{E:1}. So $(T)_{K(H)} \subsetneq (|T|)_{K(H)}$.
To see that $|T| \notin (T)_{K(H)}$, assume otherwise.
Again by Proposition \ref{P:3},
\begin{equation*}
|T| = \beta T + AT + TB + \sum_{i=1}^{m}A_{i}TB_{i}~\text{where}~ A, B, A_{i}, B_{i} \in K(H)
\end{equation*}
and so $|T|- \beta T \in (T)K(H) = (|T|)K(H)$, i.e., $\diag \left<\frac{|1-(i)^{n}\beta|}{n}\right> \in (|T|)K(H)$ and hence $\diag \left<\frac{|\beta|}{2n-1}\right> \in (|T|)K(H)$, once again contradicting Example \ref{E:1}.\\
\eE

Unlike $B(H)$-ideals, principal $K(H)$-ideals are not necessarily closed under the adjoint operation.

\bE{E: J ne J^*}
$T^{*}\notin (T)_{K(H)}$ where $T = \diag \left<\frac{(i)^n}{n}\right>$ and $T^{*} = \diag \left<\frac{(-i)^n}{n}\right>$.
Indeed, if $T^{*} \in (T)_{K(H)}$, then $\diag \left<\frac{(-i)^{n}-\alpha(i)^{n}}{n}\right> \in (T)K(H)$
implying\\
$\diag \left<\frac{|(-1)^{n}- \alpha|}{n}\right> \in  (\diag \left<\frac{1}{n}\right>)K(H)$, in particular,
$\diag \left<\frac{1}{2n}\right> \in  (\diag \left<\frac{1}{n}\right>)K(H)$.
That is, $\diag \left<\frac{1}{n}\right> \in (\diag \left<\frac{1}{n}\right>)K(H)$, contradicting Example \ref{E:1}.
\eE
\section{LATTICE STRUCTURE OF SUBIDEALS}\label{S:5}

Basic knowledge of the $B(H)$-ideal lattice contributes important perspective to operator theory.
In the spirit of recent work on this lattice by the second author joint with V. Kaftal \cite{KW11},
this section introduces the study of subideal lattices, that is, $J$-ideal lattices.

The explicit descriptions of the three types of principal $J$-ideals generated by $S$ given in Corollary \ref{C:1} implies that
\emph{every principal $J$-ideal contains $J(S)J$}.
It is well-known that every proper $B(H)$-ideal contains $F(H)$, the $B(H)$-ideal of all finite rank operators \cite[Chapter III, Section 1, Theorem 1.1]{GK69}.
So, \emph{every  nonzero principal $J$-ideal contains $F(H)$} and hence also every nonzero $J$-ideal.\\

As for the smallest proper $J$-ideal, principal or otherwise, we have $F(H)$ and it is principal.

\bR{R:6.2}
$(S)_J = F(H)$ if and only if $S \in F(H)$.
Indeed, by Corollary \ref{C:1}, $(S)_J = \mathbb C S + SJ + JS + J(S)J$.
So if $S \in F(H)$, then $(S)_J \subseteq F(H)$, and since every nonzero $J$-ideal contains $F(H)$, hence $(S)_J = F(H)$.
Conversely, if $S \in F(H)$, then $(S)_J = \mathbb C S + SJ + JS + J(S)J \subseteq F(H)$, hence $(S)_J = F(H)$.
\eR

For maximal proper $J$-ideals inside principal $J$-ideals we have:

\bP{P:MaxJidealsinsideprincipals}
Every nonzero principal $J$-ideal (all three types) generated by an operator $S$ of infinite rank contains a maximal $J$-ideal.
(See Section \ref{Questions}, Question 3).
\eP
\bp
We prove this only for principal linear $J$-ideal $(S)_J$ because the same argument holds for principal real linear $J$-ideal
$(S)_J^{\mathbb{R}}$ and principal $J$-ideal $\left<S\right>_J$.
Consider $ \mathscr A := \{J\text{-ideals}~ \mathcal I \mid \mathcal I \subsetneq (S)_J\}$.
By Corollary \ref{C:1}, principal linear $J$-ideal $(S)_J$ contains the $B(H)$-ideal $J(S)J$ which further contains $F(H)$.
Since $S$ is of infinite rank, $F(H) \subsetneq (S)_J$, and as every $B(H)$-ideal is a $J$-ideal, $F(H) \in \mathscr{A}$ and so $\mathscr A \neq \emptyset$.
Taking the inclusion partial order on $\mathscr A$ and any chain $\mathscr C$ in $\mathscr A$,
it is easy to show that $\displaystyle{\bigcup_{\mathcal I \in \mathscr C}}\mathcal I$ is a $J$-ideal contained in $(S)_J$.
It is a $J$-ideal properly contained in $(S)_J$ because it does not contain $S$.
Indeed if it did contain $S$, then there is some $S \in  \mathcal I_{o} \in \mathscr C$ implying $(S)_J \subset \mathcal I_{o}$, contradicting the criterion for  $\mathcal I_{o} \in \mathscr A$.
Therefore, every chain has an upper bound in $\mathscr A$ and so by Zorn's Lemma, $\mathscr A$ has a maximal element, i.e., $(S)_J$ has a maximal $J$-ideal.
\ep

\bR{R:7.1}
The $J$-ideal $JS + SJ + J(S)J$ is always maximal in the principal linear $J$-ideal $(S)_J$ and the principal real linear $J$-ideal $(S)_J^{\mathbb R}$,
but is never a maximal $J$-ideal in the principal $J$-ideal $\left<S\right>_J$.

Indeed, the quotient rings
$
\frac{(S)_J}{JS + SJ + J(S)J} \quad \text{and}\quad \frac{(S)_J^{\mathbb R}}{JS + SJ + J(S)J}
$
have no proper ideals because of the following reason:
Every nonzero element is a coset of the form
$\left[\alpha S\right] \in \frac{(S)_J}{JS + SJ + J(S)J}$ and $\left[rS\right] \in \frac{(S)_J^{\mathbb R}}{JS + SJ + J(S)J}$ for $\alpha \in \mathbb C,\, r \in \mathbb R$.
Therefore, the linear $J$-ideal generated by any single element $[\alpha S]$ ($\alpha \ne 0$) in the quotient ring is the entire quotient ring $\frac{(S)_J}{JS + SJ + J(S)J}$ and likewise the real linear $J$-ideal generated by any $[rS]$ is the entire quotient ring  $\frac{(S)_J^{\mathbb R}}{JS + SJ + J(S)J}$, that is, these quotient rings have no proper ideals.
From general ring theory recall that, for a fixed ideal $I$ in a ring $R$, the map $J \rightarrow [J] := \{[A] \mid A \in J\}$ is a ring inclusion preserving isomorphism (so a one-to-one correspondence) between the ideals $J$ containing $I$ and the ideals of the quotient ring $R/I$.
So because there are no proper linear $J$-ideals in the quotient ring $\frac{(S)_J}{JS + SJ + J(S)J}$ and real linear $J$-ideals in
$\frac{(S)_J^{\mathbb R}}{JS + SJ + J(S)J}$,
there are no proper $J$-ideals (neither linear inside $(S)_J$ nor real linear inside $(S)_J^{\mathbb R}$) containing $JS + SJ + J(S)J$.
That is, $JS + SJ + J(S)J$ is always a maximal linear $J$-ideal in $(S)_J$ and a maximal real linear $J$-ideal in $(S)_J^{\mathbb R}$.

But in the quotient ring $\frac{\left<S\right>_J}{JS + SJ + J(S)J}$ every non-zero element is of the form $[nS]$ for $n \in \mathbb Z \setminus \{0\}$,
and when $n \in \mathbb Z \setminus \{0, \pm 1\}$ the ideal $\mathcal J$ generated by $[nS]$ in the quotient ring is proper because the class of integer multiples of $[nS]$ never contains $[S]$.
Therefore the inverse image of $\mathcal J$ under the natural projection map
$\pi :\left<S\right>_J \rightarrow \frac{\left<S\right>_J}{JS + SJ + J(S)J}$ is a proper ideal in  $\left<S\right>_J$ containing $JS + SJ + J(S)J$.
Therefore, $JS + SJ + J(S)J$ is never a maximal $J$-ideal in $\left<S\right>_J$.
\eR

Considering Remark \ref{R:6.2} and Remark \ref{R:7.1},
we wonder if $S$ being finite rank is necessary for $JS + SJ + J(S)J$ to be principal? (See Section \ref{Questions}, Question 5.)

\bR{R:7.2}
For idempotent $B(H)$-ideals $J$ ($J^{2} = J$), since $SJ, JS \subseteq J(S)$, one has $J(S)J = SJ + JS + J(S)J$.
We wonder if being idempotent is necessary for this equality? (See Section \ref{Questions}, Question 6.)
\eR

As mentioned earlier, finitely generated $B(H)$-ideals are principal, but this is not generally the case as seen in the next proposition.
This justifies Example \ref{E:4} and its preceding comment.
\bP{P:7.2}
If $S, T$ are simultaneously diagonalizable operators with disjoint supports and $s(S) \cong s(T)$ ($as(S) \leq  s(T) \leq bs(S)$ for some $a,b > 0$),
then the linear $J$-ideal $(\{S,T\})_J$ is a principal linear $J$-ideal if and only if $(\{S,T\})$ is $J$-soft.
\eP
\bp
($\Rightarrow$):
If $(\{S,T\})_J = (A)_J$ for some $A \in J$, then $S, T \in (A)_J$ so, by Corollary \ref{C:1} it is immediate that $S = \alpha A + X$ and $T = \beta A + Y$ for some $X, Y \in J(A),~ \alpha, \beta \in \mathbb C$.
And since $A \in (\{S,T\})_J$, by Proposition \ref{P:4}, one has $A = c_{1}S + c_{2}T + R$ for some $R \in  J(\{S,T\}),~c_{1}, c_{2} \in \mathbb C$.
If $\alpha = \beta = 0$, then $S, T \in J(A) = J( c_{1}S + c_{2}T + R) \subseteq J(\{S,T\})$ implying $(\{S,T\}) \subseteq J(\{S,T\})$,
and since the reverse inclusion is automatic, one obtains $J$-softness of $(\{S,T\})$.
Otherwise, $-\beta S + \alpha T = -\beta X + \alpha Y \in J(\{S,T\})$.
Since $S$ and $T$ are simultaneously diagonalizable with diagonals of disjoint support,
there is a unitary operator $U$ such that $U^{*}SU = D_{1}$ and $U^{*}TU = D_{2}$ where $D_{1}$ and $D_{2}$ have disjoint supports.
Since $J(\{S,T\})$ is a $B(H)$-ideal, $U^{*}(-\beta S + \alpha T)U  = -\beta D_{1} + \alpha D_{2} \in J(\{S,T\})$.\\
If $\alpha = \beta = 0$, then $S = X, T = Y \in J(A) \subseteq J(\{S,T\})$, i.e., $(\{S,T\}) \subseteq J(\{S,T\})$ so $(\{S,T\})$ is $J$-soft.
If $\alpha = 0, \beta \ne 0$, then $S = X \in J(A) \subseteq J(\{S,T\})$ and $-\beta D_{1} \in J(\{S,T\})$ and hence $S \in J(\{S,T\})$.
Since $s(T) \leq bs(S)$, the hereditary property of $\Sigma (J(\{S,T\}))$ implies $\diag s(T) \in J(\{S,T\})$
and hence $T \in J(\{S,T\})$ so $(\{S,T\})$ is $J$-soft.
Similarly if $\beta = 0, \alpha \ne 0$.
Finally if both $\alpha, \beta \ne 0$, then as $D_{1}, D_{2}$ have disjoint supports,
$|\beta s(D_{1})|, |\alpha s(D_{2})| \leq s(-\beta D_{1} + \alpha D_{2})$
implying $D_{1}, D_{2} \in J(\{S,T\})$ and then $S, T \in J(\{S,T\})$, hence $(\{S,T\}) = J(\{S,T\})$.

$\Leftarrow$: By Theorem \ref{T:3}, $J$-softness of $(\{S,T\})$ implies $(\{S,T\})_J$ is a $B(H)$-ideal and Lemma \ref{L:3.2} implies $(\{S,T\})_J = (|S| + |T|)$ which is a principal $B(H)$-ideal.
But because $(\{S,T\})_J = (|S| + |T|)$, the principal $B(H)$-ideal $(|S| + |T|)$ is $J$-soft.
Hence $(\{S,T\})_J = (|S| + |T|)_J$ since
$(|S| + |T|)_J = (|S| + |T|)$ by Theorem \ref{T:2}.
\ep

\section{QUESTIONS}\label{Questions}

Natural questions arise from this work.\\

\noindent 1. Is it true that if $\mathscr S$ is countable, then $(\mathscr S)_J$ is a $B(H)$-ideal if and only if $(\mathscr S)$ is $J$-soft,
i.e., $(\mathscr S) = J(\mathscr S)$? (See Sections 3-4.)\\

\noindent 2. Find necessary and sufficient conditions for when a finitely generated $J$-ideal is a principal $J$-ideal? \\

\noindent Proposition \ref{P:7.2} indicates that $J$-softness is a necessary and sufficient condition for the particular class:
$S, T$ are mutually diagonalizable operators with disjoint supports and equivalent s-number sequences.\\

\noindent 3. Are maximal ideals inside principal $J$-ideals (in particular, those guaranteed by Proposition \ref{P:MaxJidealsinsideprincipals})
always principal or always non-principal? \\

\noindent 4. Find necessary and sufficient conditions for two principal $J$-ideals to be equal?
(Examples 5.1-5.2 help motivate this natural but deceptively simple question.) \\

\noindent 5. Can $JS + SJ + J(S)J$ be a principal $J$-ideal besides the case $S \in F(H)$? \\
(See Remarks \ref{R:6.2} and \ref{R:7.1}.)\\

\noindent 6. Find a necessary and sufficient condition(s) to make $J(S)J = JS + SJ + J(S)J$. (Recall Remark \ref{R:7.2}.)\\

\noindent 7. Does $(\mathscr S)_J^\mathbb R = (\mathscr S)_J$ imply $(\mathscr S)$ is $J$-soft?\\
(See comment succeeding Theorem \ref{T:4}).\\

\noindent \emph{Acknowledgement.} The second author thanks David Pitts for input on this subject during a visit in 1998.

\end{document}